\documentclass[12pt]{article}

\usepackage{amssymb}
\usepackage{amsmath,amsthm}
\usepackage[dvips]{graphicx}
\usepackage{float}

\newtheorem{thm}{Theorem}[section]
\newtheorem{cor}[thm]{Corollary}
\newtheorem{lem}[thm]{Lemma}

\theoremstyle{definition}

\newtheorem{eg}[thm]{Example}
\theoremstyle{remark}

\newtheorem*{ack}{Acknowledgments}

\title{{\Large The complete splitting number of a lassoed link}}
\author{Ayaka Shimizu \\ \footnotesize{Graduate School of Science, Osaka City University,} \\[-3pt]
\footnotesize{Sugimoto, Sumiyoshi-ku Osaka 558-8585, Japan}\\
\footnotesize{Email: shimizu1984@gmail.com}}

\begin{document}

\maketitle

\begin{abstract}
In this paper, we define a lassoing on a link, a local addition of a trivial knot to a link. 
Let $K$ be an $s$-component link with the Conway polynomial non-zero. 
Let $L$ be a link which is obtained from $K$ by $r$-iterated lassoings. 
The complete splitting number $\mathrm{split}(L)$ is greater than or equal to $r+s-1$, 
and less than or equal to $r+\mathrm{split}(K)$. 
In particular, we obtain from a knot by $r$-iterated component-lassoings an algebraically completely splittable link $L$ with $\mathrm{split}(L)=r$. 
Moreover, we construct a link $L$ whose unlinking number is greater than $\mathrm{split}(L)$. 
\end{abstract}

\section{Introduction}

\noindent The splittability of a link is one of the basic concepts in knot theory. 
For example, the splittability interacts with polynomial invariants: 
the Alexander polynomial and the Conway polynomial take zero for a splittable link. 
Jones polynomial and skein polynomial have specific formulae with respect to the split sum. 
Moreover, the splittabilities of links or spatial graphs are studied and applied to other subjects: chemistry, biology, psychology, etc. 
For example, Kawauchi proposed a model of prion proteins as a spatial graph \cite{kawauchip2}, and 
Yoshida studied its splittability which concerns with the study of prion diseases: mad cow disease, scrapie, Creutzfeldt-Jakob disease, etc. \cite{yoshida}. 
Another example is about a model of human mind which is also proposed by Kawauchi \cite{kawauchip}, \cite{kawauchip2}; 
by considering one's mind as a knot and by considering a mind relation of $n$ persons as an $n$-component link, 
the models ``mind knots'' and ``mind links'' are studied. 
The splittability of a link corresponds to the ``self-releasability'' of a mind link. 

\begin{figure}[!ht]
\begin{center}
\includegraphics[width=50mm]{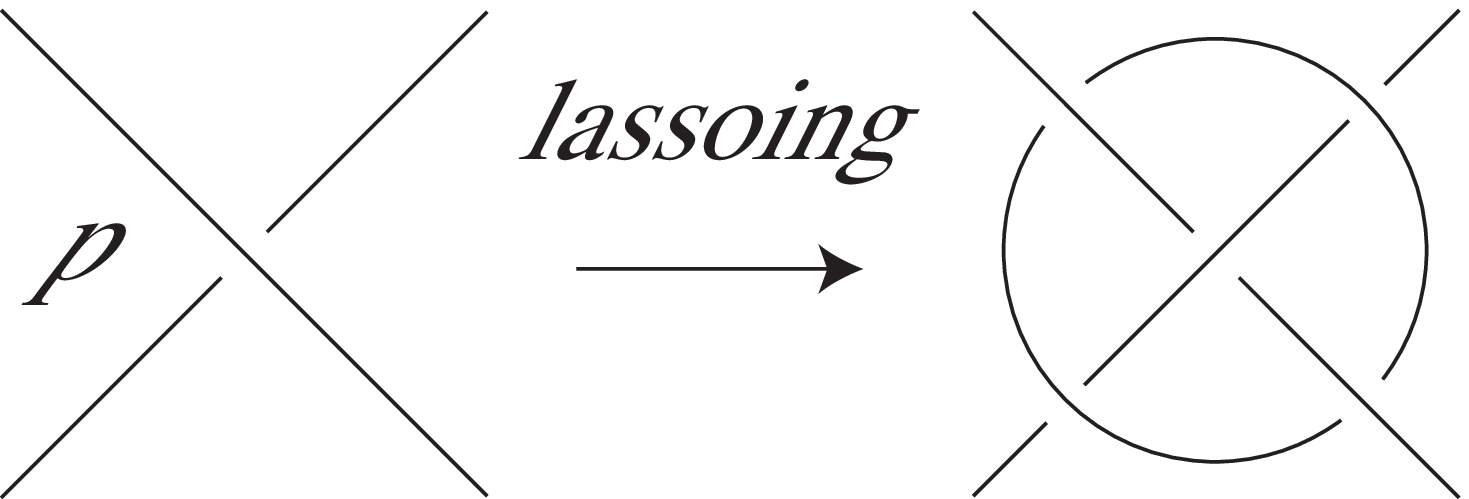}
\caption{}
\label{quoitting}
\end{center}
\end{figure}

For a two-component link, Adams defined the splitting number which represents how distant the link is from a splittable link \cite{adams}. 
In this paper, we define for an $n$-component link $L$ $(n=2,3,4,\dots )$ the complete splitting number $\mathrm{split}(L)$ which represents how distant the link is from a completely splittable link. 
The \textit{unlinking number} $u(L)$ of a link $L$ is the minimal number of crossing changes in any diagram of $L$ which are needed to obtain the trivial link $L$. 
Since a trivial link is completely splittable, we have $\mathrm{split}(L)\le u(L)$. 
Lassoing is a crossing-changing and loop-adding local move as shown in Figure \ref{quoitting} 
(we give the precise definitions of completely splittable, complete splitting number, and a lassoing in Section \ref{csn}). 

\noindent For any $r$-component link $L=L_1\cup L_2\cup...\cup L_r$ $(r=1,2,3,\dots )$ with the Conway polynomial 
$\nabla(L)\ne 0$, there are $\ell$-iterated lassoings from $L$ to an
algebraically completely splittable link $L^*$ with $\nabla (L^*)\ne 0$ where
$\ell =\sum_{i<j}|Link(L_i,L_j)|$ (we define an algebraically completely splittable link in Section \ref{csn}). 
For any $s$-component link $K=K_1\cup K_2\cup...\cup K_s$ $(s\ge 1)$ with 
$\nabla(K)\ne 0$, there are $(\ell +u)$-iterated lassoings from $K$ to an
algebraically completely splittable link L with trivial components such that
$\nabla(L)\ne 0$ where
$\ell =\sum_{i<j}|Link(K_i,K_j)|$ and $u=\sum_{i=1}^s u(K_i)$.
In this paper, we show the following theorem: 

\phantom{x}
\begin{thm}
Any link $L$ obtained from any $s$-component link $K=K_1\cup K_2\cup...\cup K_s$ $(s=1,2,3,\dots )$ with
$\nabla(K)\ne 0$ by $r$-iterated lassoings $(r=0,1,2,\dots )$ satisfies 
$$r+split(K)\ge split(L)\ge r+s-1.$$
\label{mainthm}
\end{thm}
\phantom{x}

\noindent We have the following corollaries: 

\phantom{x}
\begin{cor}
For any $s$-component link $K=K_1\cup K_2\cup...\cup K_s$ $(s=1,2,3,\dots )$ with $\mathrm{split}(K)=s-1$, 
and any integer $r\ge \ell +u$ where $\ell =\sum_{i<j}|Link(K_i,K_j)|$ and $u=\sum_{i=1}^s u(K_i)$, 
there are $r$-iterated lassoings from $K$ to an algebraically completely splittable link $L$ with trivial components such that 
$\mathrm{split}(L)=r+s-1$. 
\end{cor}
\phantom{x}

\phantom{x}
\begin{cor}
Let $K$ be a knot. Let $L$ be a link which is obtained from $K$ by $r$-iterated lassoings ($r=1,2,3,\dots $). 
Then $L$ has $\mathrm{split}(L)=r$. 
\label{knot-lasso}
\end{cor}
\phantom{x}

\noindent We define a \textit{component-lassoing} to be the lassoing at a self-crossing point of a diagram. 
We have the following corollary: 

\phantom{x}
\begin{cor}
Every link $L$ obtained from a knot $K$ by $r$-iterated component-lassoings $(r=1,2,3,\dots )$ 
is an $(r+1)-$component algebraically completely splittable link with $\mathrm{split}(L)=r$. 
\label{comp-lasso-cor}
\end{cor}
\phantom{x}

\begin{figure}[!ht]
\begin{center}
\includegraphics[width=30mm]{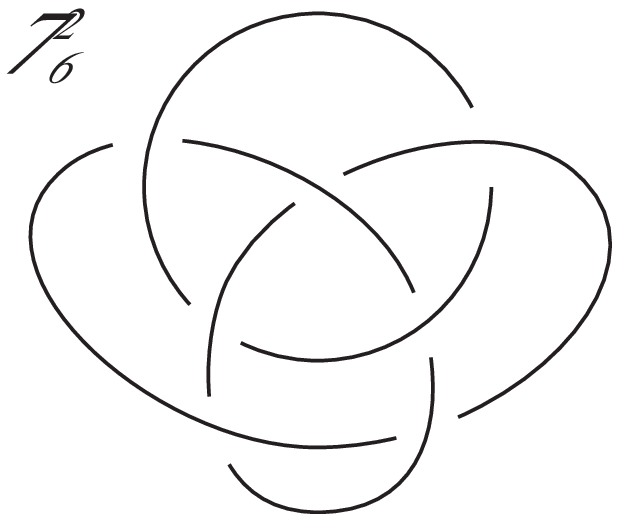}
\caption{}
\label{7-2-6}
\end{center}
\end{figure}

\noindent For example, the link $7^2_6$ depicted in Figure \ref{7-2-6} which is a link obtained from a trefoil knot by a single component-lassoing, 
has the linking number zero and $\mathrm{split}(7^2_6)=1$. 
We also remark that $u(7^2_6)=2$ (\cite{kohn}). 
Adams also showed in \cite{adams} that there is a two-component link, each component of which is trivial, but such that its splitting number is less than its unlinking number, like the link $7^2_6$. 
We show in Section \ref{n-t} that for any integer $r>0$ and any knot $K$ with Nakanishi's index $e(K)>2r$, 
any link $L$ obtained from $K$ by $r$-iterated lassoings is a link such that $\mathrm{split}(L)<u(L)$, i.e., $L$ is non-trivial by any $r$ crossing changes.

\section{Complete splitting number}
\label{csn}

\begin{figure}[!ht]
\begin{center}
\includegraphics[width=60mm]{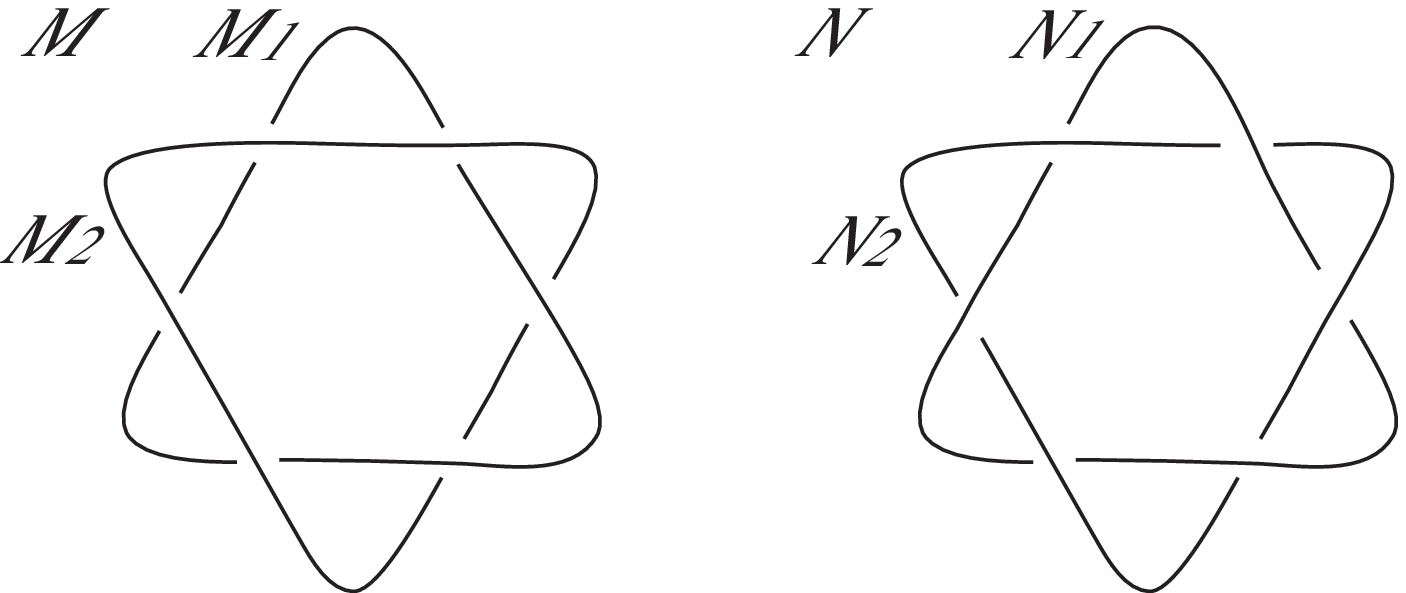}
\caption{}
\label{splitinto}
\end{center}
\end{figure}

\noindent Let $L=L_1\cup L_2\cup \dots \cup L_r$ be a link consisting of sublinks $L_i$ $(i=1,2,\dots ,r)$. 
A link $L$ is \textit{splittable into} $L_1,L_2,\dots ,L_r$ if there exist mutually disjoint 3-balls $B_i$ ($i=1,2,\dots ,r$) 
in $S^3$ such that $L_i \subset B_i$. 
For example, the link $M$ in Figure \ref{splitinto} is splittable into $M_1$ and $M_2$ 
whereas the link $N$ is not splittable into $N_1$ and $N_2$. 
A link $L$ is \textit{splittable} if $L$ is splittable into subdiagrams $L_1$ and $L_2$, where $L=L_1\cup L_2$, $L_1,L_2\neq \phi$. 
For example, the link $M$ in Figure \ref{splitinto} is a splittable link. 

\begin{figure}[!ht]
\begin{center}
\includegraphics[width=60mm]{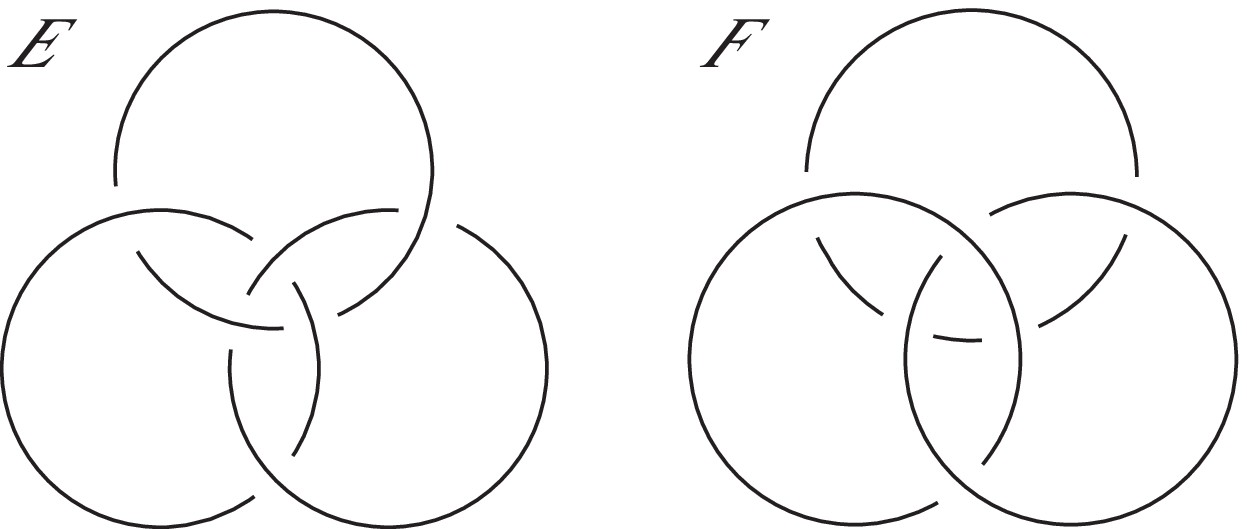}
\caption{}
\label{compsplit}
\end{center}
\end{figure}

\noindent A link $L$ is \textit{completely splittable} if $L$ is splittable into all the knot components of $L$. 
In particular, a knot is assumed as a non-splittable link but a completely splittable link. 
A link $L$ is \textit{algebraically completely splittable} if every two knot components $K_i$ and $K_j$ of $L$ have the linking number 
$\mathrm{Link}(K_i,K_j)=0$. 
For example, the link $E$ in the left hand of Figure \ref{compsplit} is not completely splittable 
but algebraically it is completely splittable. 

\noindent The \textit{complete splitting number} $\mathrm{split}(D)$ of a link diagram $D$ is the minimal number of crossing changes 
which are needed to obtain a diagram of a completely splittable link from $D$. 
For example, the link diagram $F$ in the right hand in Figure \ref{compsplit} has $\mathrm{split}(F)=1$. 
As a relation to the warp-linking degree $ld(D)$ of $D$, we have $\mathrm{split}(D)\leq ld(D)$, 
where the warp-linking degree is a restricted warping degree which can be calculated directly or by using matrices \cite{shimizu-d}, \cite{shimizu}. 
The \textit{complete splitting number} $\mathrm{split}(L)$ of a link $L$ is the minimal number of crossing changes in any diagram of the link which are needed to obtain a completely splittable link.

\begin{figure}[!ht]
\begin{center}
\includegraphics[width=80mm]{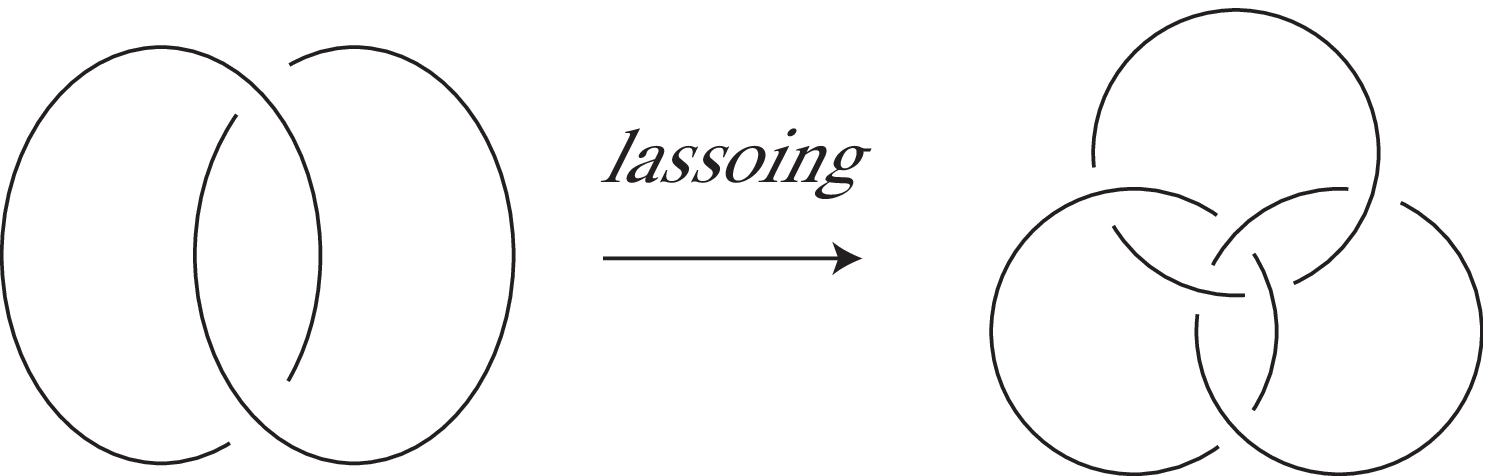}
\caption{}
\label{lasso-ex}
\end{center}
\end{figure}

\noindent Let $p$ be a crossing point of a link diagram $D$. 
We put a lasso around $p$, i.e., we apply a crossing change at $p$, and add a loop alternately around the crossing 
as shown in Figure \ref{quoitting}. 
Then, we obtain another link diagram $D'$. 
The diagram $D'$ is said to be obtained from $D$ by \textit{lassoing} at $p$. 
Let $L'$ be the link which has the diagram $D'$. 
The link $L'$ is said to be obtained from $L$ by a \textit{lassoing}. 
For example, we obtain the Borromean ring from the Hopf link by a lassoing (see Figure \ref{lasso-ex}). 

\begin{figure}[!ht]
\begin{center}
\includegraphics[width=145mm]{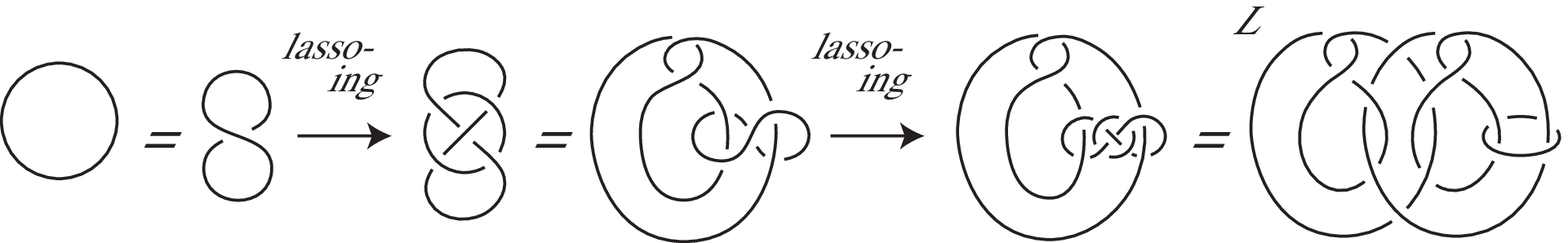}
\caption{}
\label{2-lasso}
\end{center}
\end{figure}

\noindent A link $L'$ is said to be obtained from $L$ by \textit{$r$-iterated lassoings} 
if $L'$ is obtained from $L$ by lassoings $r$ times iteratively. 
For example, the link $L$ in Figure \ref{2-lasso} is a link obtained from a trivial knot by two-iterated lassoings. 
Since a lassoing depends on the choice of a crossing point and the choice of a diagram of the link, 
we may have many types of link by a lassoing.

\section{Conway polynomial}

\noindent Let $\nabla (L;z)$ be the Conway polynomial of a link $L$ with an orientation. 
We have the following lemma: 

\phantom{x}
\begin{lem}
We have 
\begin{align*}
\nabla \left( \begin{minipage}{20pt}
              \includegraphics[width=20pt]{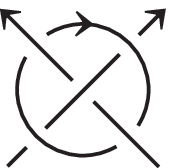}
         \end{minipage}
;z\right) =-z^3\nabla \left( \begin{minipage}{20pt}
              \includegraphics[width=18pt]{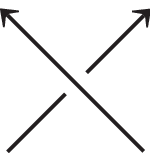}
         \end{minipage}
;z\right) , 
\end{align*}

\begin{align*}
\nabla \left( \begin{minipage}{20pt}
              \includegraphics[width=20pt]{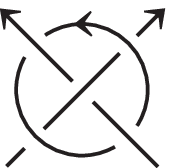}
         \end{minipage}
;z\right) =-z^3\nabla \left( \begin{minipage}{20pt}
              \includegraphics[width=18pt]{12.eps}
         \end{minipage}
;z\right) , 
\end{align*}

\begin{align*}
\nabla \left( \begin{minipage}{20pt}
              \includegraphics[width=20pt]{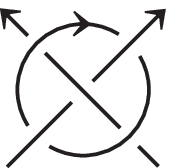}
         \end{minipage}
;z\right) =z^3\nabla \left( \begin{minipage}{20pt}
              \includegraphics[width=18pt]{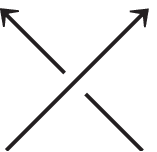}
         \end{minipage}
;z\right) ,
\end{align*}

\begin{align*}
\nabla \left( \begin{minipage}{20pt}
              \includegraphics[width=20pt]{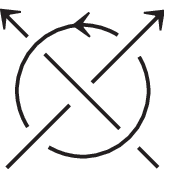}
         \end{minipage}
;z\right) =z^3\nabla \left( \begin{minipage}{20pt}
              \includegraphics[width=18pt]{32.eps}
         \end{minipage}
;z\right) . 
\end{align*}
\phantom{x}
\begin{proof} 
We obtain the first equality by the skein relations in Figure \ref{conway}. 

\begin{figure}[!ht]
\begin{center}
\includegraphics[width=115mm]{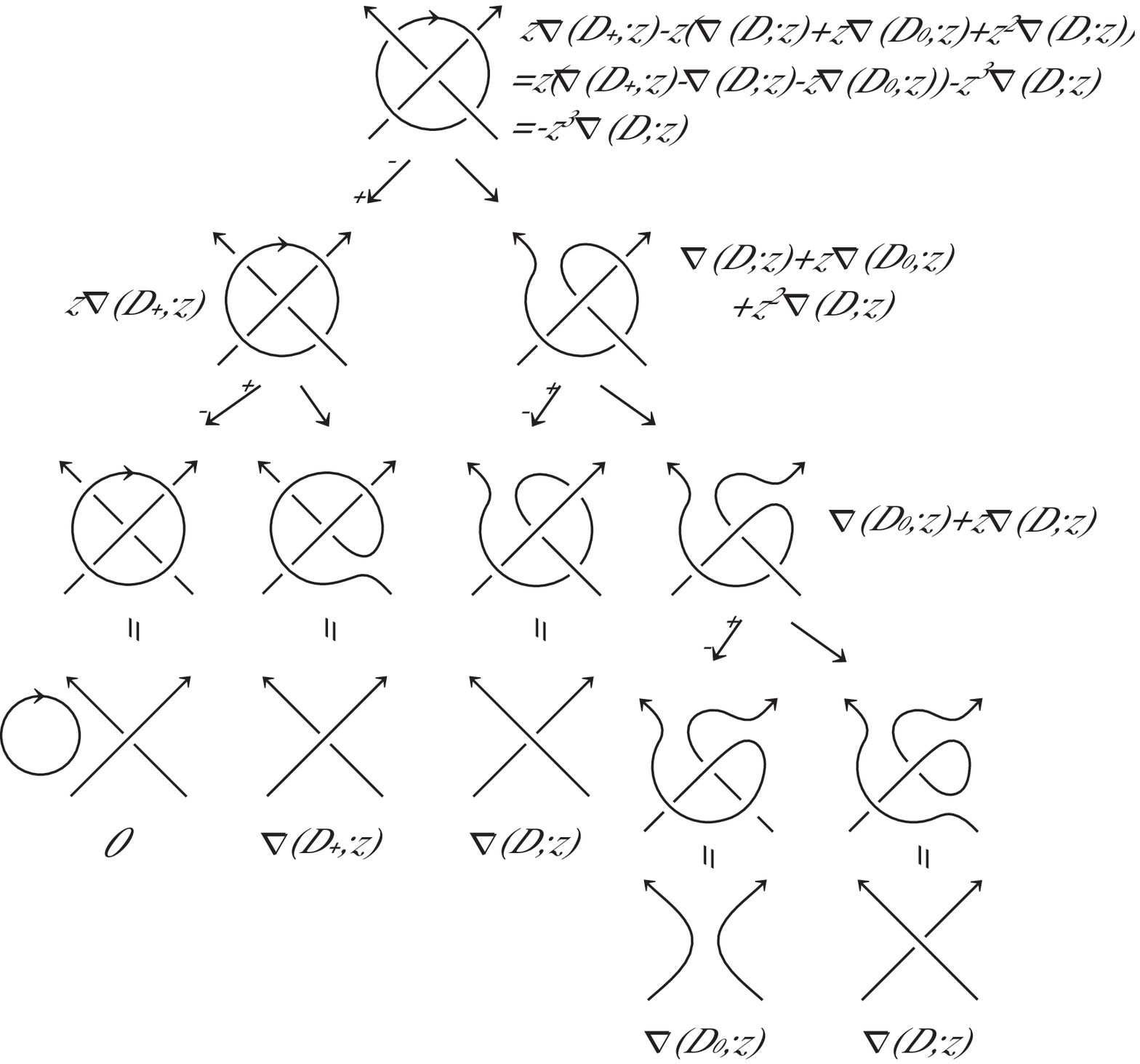}
\caption{}
\label{conway}
\end{center}
\end{figure}

\noindent The other equalities are similarly obtained.
\end{proof}
\label{lem-conway}
\end{lem}
\phantom{x}

\begin{eg}
The link diagram $D$ in Figure \ref{conway-ex} is obtained from a diagram of a trefoil knot by 2-iterated lassoings. 
Then we have $\nabla (L)=z^3\times z^3\times \nabla (3_1)=z^6(1+z^2)$, where $L$ is a link represented by $D$, and $3_1$ is a trefoil knot. 

\begin{figure}[!ht]
\begin{center}
\includegraphics[width=25mm]{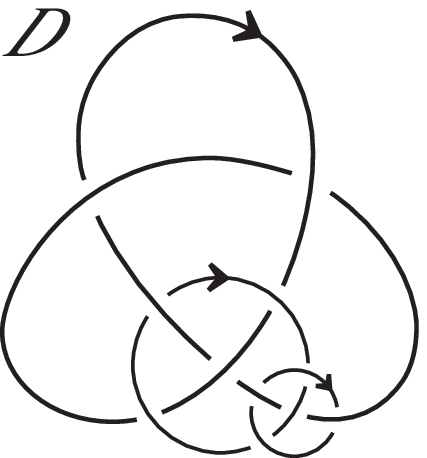}
\caption{}
\label{conway-ex}
\end{center}
\end{figure}
\end{eg}
\phantom{x}

\noindent We remark that for a link $L'$ with $\nabla (L')=0$, there are no lassoings from $L'$ to $L$ with $\nabla (L)\neq 0$. 
We have the following corollary: 

\phantom{x}
\begin{cor}
Let $L$ be a link obtained from a link $L'$ with $\nabla (L')\neq 0$, in particular from any knot $K$, 
by $r$-iterated lassoings ($r=1,2,3,\dots $). 
Then we have $\nabla (L)\neq 0$. 
\label{knot-quoit}
\end{cor}
\phantom{x}

\noindent Let $\Lambda$ be the integral Laurent polynomial ring, i.e., $\Lambda =\mathbb{Z}[t,t^{-1}]$. 
With respect to the one-variable Alexander polynomial, we have the following corollary by Lemma \ref{lem-conway} 
by substituting $t^{\frac{1}{2}}-t^{-\frac{1}{2}}$ for $z$:

\phantom{x}
\begin{cor}
Let $L'$ be a link which is obtained from a link $L$ by a lassoing. Then we have 
$$\Delta _{L'}(t){\dot =}(t-1)^3\Delta _L(t), $$
where $\Delta _L(t)$ is the one-variable Alexander polynomial of $L$, 
and ${\dot =}$ means equal up to multiplications of the units of $\Lambda$. 
\end{cor}
\phantom{x}

\noindent We show an example. 

\phantom{x}
\begin{eg}
We have 
$$\Delta (5^2_1)\doteq \Delta (7^2_8)\doteq \Delta (8^2_{15})\doteq (t-1)^3\Delta (0),$$
$$\Delta (7^2_6)\doteq \Delta (9^2_{55})\doteq \Delta (9^2_{56})\doteq (t-1)^3\Delta (3_1),$$
$$\Delta (8^2_{13})\doteq (t-1)^3\Delta (4_1),$$
$$\Delta (9^2_{31})\doteq (t-1)^3\Delta (5_1),$$
$$\Delta (9^2_{32})\doteq \Delta (9^2_{33})\doteq (t-1)^3\Delta (5_2),$$
where $\Delta (L)=\Delta _L(t)$. 
All the two-component links with the crossing number nine or less which are obtained from knots by lassoings 
have been listed above. 
\end{eg}
\phantom{x}

\noindent Up to multiplications of $t-1$, the one-variable Alexander polynomial of any link is the Alexander polynomial of an algebraically completely 
splittable link consisting of trivial components: 

\phantom{x} 

\begin{cor}
Let $(t-1)^mf(t)$ be the Alexander polynomial of a link, where $m$ is a non-negative integer, $f(t)\in \Lambda$, and $f(1)\neq 0$. 
Then, there exists a non-negative integer $n$ such that the Laurent polynomial $(t-1)^{m+3n}f(t)$ is the Alexander polynomial of 
an algebraically completely splittable link consisting of trivial components. 

\begin{proof}
We can change a crossing by a lassoing.
\end{proof}
\end{cor}
\phantom{x}

\section{Proof of Theorem \ref{mainthm}}

\noindent In this section, we prove Theorem \ref{mainthm}. 
Before the proof, we define some notions which are due to \cite{kawauchi} to prove Theorem \ref{mainthm}. 
For the integral Laurent polynomial ring $\Lambda =\mathbb{Z}[t,t^{-1}]$, 
a \textit{multiplicative set} of $\Lambda$ is a subset $S\subset \Lambda -\{ 0\} $ which satisfies the following three conditions: 
the units $\pm t^i$ ($i\in \mathbb{Z}$) are in $S$, 
the product $gg'$ of any elements $g$ and $g'$ of $S$ is in $S$, and 
every prime factor of any element $g\in S$ is in $S$. 
For the quotient field $Q(\Lambda)$ of $\Lambda$ and a multiplicative set $S$ of $\Lambda$, 
$\Lambda _S=\{ f/g \in Q(\Lambda ) | f \in \Lambda , g \in S \}$ is a subring of $Q(\Lambda)$. 
For a finitely generated $\Lambda$-module $H$, let $H_S$ be the $\Lambda _S$-module $H\otimes _{\Lambda}\Lambda _S$, and 
$e_S(H)$ the minimal number of $\Lambda _S$-generators of $H_S$. 
We take $e_S(H)=0$ when $H=0$. 
We call $e_{S}(H)$ the $\Lambda _{S}$-\textit{rank} of $H$. 
Let $L$ be an oriented link in $S^3$, and $E(L)=cl(S^3-L)$ the compact exterior of $L$. 
Let $\tilde{E}(L)\rightarrow E(L)$ be the infinite cyclic covering 
which is induced from the epimorphism $\gamma _L:\pi _1(E(L))\rightarrow \mathbb{Z}$ 
sending each oriented meridian of $L$ to $1\in \mathbb{Z}$.  
Then we can regard $H_1(\tilde{E}(L))$ as a finitely $\Lambda$-module. 
We denote $e_S(H_1(\tilde{E}(L)))$ by $e_S(L)$. 
Let $L$, $L'$ be links which have the same number of components. 
By Theorem 2.3 in \cite{kawauchi}, we immediately have 
\begin{align}
d^X(L,L')\geq \vert e_s(L)-e_s(L')\vert ,
\label{ka23}
\end{align}
where $d^X(L,L')$ denotes the $X$-distance between $L$ and $L'$. 
We prove Theorem \ref{mainthm}. 

\phantom{x}
\noindent {\it Proof of Theorem \ref{mainthm}.}
Let $L$ be a link which is obtained from a link $K=K_1\cup K_2\cup \dots \cup K_s$ with $\nabla (K)\ne 0$ 
by $r$-iterated lassoings $(r=1,2,3,\dots )$. 
Let $L'$ be a completely splittable link which is obtained from $L$ by $m$ crossing changes, 
where $m=\mathrm{split}(L)=d^X(L,L')$. 
We set $S=\Lambda -\{0\} $. 
Since $L'$ is completely splittable and the number of components of $L'$ is $r+s$, we have 
\begin{align}
e_s(L')=r+s-1. 
\label{esl0}
\end{align}
The Alexander polynomial of $L$ is non-zero because the Conway polynomial of $L$ is non-zero by Corollary \ref{knot-quoit}. 
Hence we have 
\begin{align}
e_s(L)=0. 
\label{esl1}
\end{align}
By substituting the equalities (\ref{esl0}), (\ref{esl1}) and $d^X(L,L')=\mathrm{split}(L)$ into the inequality (\ref{ka23}), we have 
$$\mathrm{split}(L)\geq r+s-1.$$
From the $r$-iterated lassoings, we have 
$$r+\mathrm{split}(K)\ge \mathrm{split}(L).$$
Hence we have the inequality
$$r+\mathrm{split}(K)\ge \mathrm{split}(L)\ge r+s-1.$$
\hfill$\square$\\
\phantom{x}

\noindent As the contraposition of Theorem \ref{mainthm}, we have the following corollary: 

\phantom{x}
\begin{cor}
Let $K=K_1\cup K_2\cup \dots \cup K_s$ be an $s$-component link. 
If $K$ has $\mathrm{split}(K)<s-1$, then $\nabla (K)=0$. 
\end{cor}

\section{Non-triviality}
\label{n-t}

\noindent In this section, we discuss the non-trivialities of completely splittable links which are obtained from $L$ 
in Corollary \ref{knot-lasso} by $r$ crossing changes $(r=1,2,\dots )$. 
For a link $L$ obtained from a knot $K$ by $r$-iterated lassoings, we have the following theorem: 

\phantom{x}
\begin{thm}
If a link $L$ is obtained from a knot $K$ with $e(K)>2r$ by $r$-iterated lassoings $(r=1,2,3,\dots )$, 
then we have $\mathrm{split}(L)=r$ and $u(L)>r$. 
\label{split-and-u}
\end{thm}
\phantom{x}

\noindent We remark that in Theorem \ref{split-and-u} the link L is an algebraically completely splittable link if the $r$-iterated lassoings are all component-lassoings. 
Before proving Theorem \ref{split-and-u}, we have the following Lemma:

\phantom{x}
\begin{lem}
Let $L_0=K_1+K_2+\dots +K_{r}$ be a completely splittable link with $r$ components. 
Then we have 
$$u(L_0)=\sum _{i=1}^ru(K_i).$$

\phantom{x}
\begin{proof}
We have $u(L_0)=u+u_1+u_2+\dots +u_r$, 
where $u$ is the number of non-self crossing changes and 
$u_i$ is the number of crossing changes on $K_i$ which are needed to obtain the trivial link from $L_0$. 
Then we have 
$$u(L_0)=u+u_1+\dots +u_r\ge u_1+\dots +u_r\ge \sum _{i=1}^ru(K_i).$$ 
Since $L_0$ is completely splittable, we have 
$$u(L_0)\le \sum _{i=1}^ru(K_i).$$
Therefore the equality holds. 
\end{proof}
\label{lem-u}
\end{lem}
\phantom{x}

\noindent We show Theorem \ref{split-and-u}

\begin{figure}[!ht]
\begin{center}
\includegraphics[width=50mm]{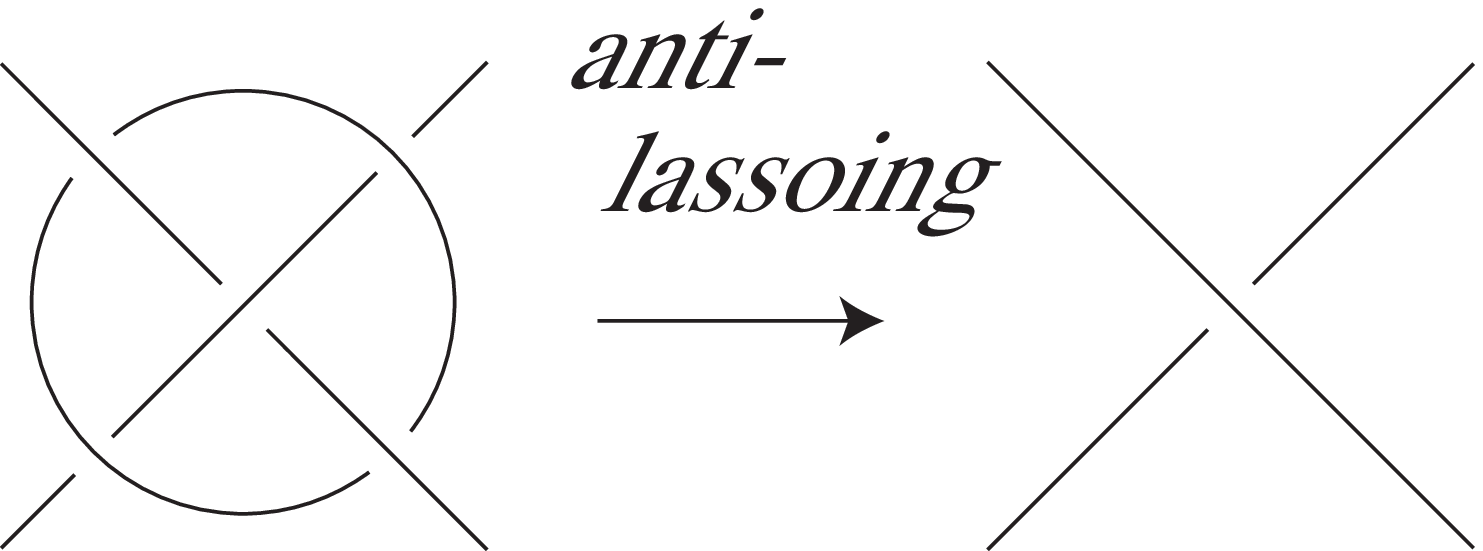}
\caption{}
\label{anti-lasso}
\end{center}
\end{figure}

\phantom{x}
\noindent {\it Proof of Theorem \ref{split-and-u}.}
Let $L_0=K_1+K_2+\dots +K_{r+1}$ be a completely splittable link 
which is obtained from $L$ by $r$ crossing changes. 
For the integral Laurent polynomial ring $\Lambda =\mathbb{Z}[t,t^{-1}]$, 
we take the multiplicative set $S$ of $\Lambda$ so that $S$ is the set of units of $\Lambda$. 
Then $e_S(L)$ is equivalent to Nakanishi's index $e(L)$ \cite{kawauchi}. 
Since we can consider $L_0=K_1+K_2+\dots +K_{r+1}$ to be a connected sum $O^{r+1}\# K_1\# K_2\# \dots \# K_{r+1}$, we have 
\begin{align*}
H_1(\tilde{E}(L_0))&\cong H_1(\tilde{E}(O^{r+1}))\oplus H_1(\tilde{E}(K_1\# K_2\# \dots \# K_{r+1})) \\
                   &\cong \Lambda ^r \oplus H_1(\tilde{E}(K_1\# K_2\# \dots \# K_{r+1})). 
\end{align*}
And by \cite{kawauchi2}, we have $e(L_0)=r+e(K_1\# K_2\# \dots \# K_{r+1})$. 
By substituting this into the inequality (\ref{ka23}), we have 
\begin{align*}
d^X(L,L_0)\geq |e(L)-e(L_0)|\geq e(L)-r-e(K_1\# K_2\# \dots \# K_{r+1}).
\end{align*}
Recall that $d^X(L,L_0)=\mathrm{split}(L)=r$. Then we have 
\begin{align}
r\geq e(L)-r-e(K_1\# K_2\# \dots \# K_{r+1}).
\label{byconnectedsum}
\end{align}
Next, we consider another completely splittable link $K+O^r$ which is obtained from $L$ by the $r$ anti-lassoings (see Figure \ref{anti-lasso}). 

\noindent Since $K+O^r=O^{r+1}\# K$, we have 
\begin{align*}
H_1(\tilde{E}(L_0))\cong \Lambda ^r \oplus H_1(\tilde{E}(K)). 
\end{align*}
And by \cite{kawauchi2}, we have $e(K+O^r)=r+e(K)$. Hence we have 
\begin{align}
r\geq r+e(K)-e(L)
\label{byantilassoings}
\end{align}
by \cite{kawauchi}. By summing the inequalities (\ref{byconnectedsum}) and (\ref{byantilassoings}), we have 
$$2r\geq e(K)-e(K_1\# K_2\# \dots \# K_{r+1}).$$
From Lemma \ref{lem-u}, we have 
$$u(L_0)=\sum _{i=1}^{r+1}u(K_i)\geq u(K_1\# K_2\# \dots \# K_{r+1})\geq e(K_1\# K_2\# \dots \# K_{r+1})\geq e(K)-2r.$$
Hence $L_0$ is non-trivial if $e(K)>2r$. 
\hfill$\square$\\
\phantom{x}

\noindent For a knot which has Nakanishi's index large enough, we can construct a link such that the unlinking number is 
greater than the complete splitting number. Here is an example. 

\phantom{x}
\begin{eg}
Since the knot $K$ in Figure \ref{connected-tref} which is the connected sum of $2r+1$ trefoil knots has Nakanishi's index $e(K)=2r+1$, 
any link $L$ obtained from $K$ by $r$-iterated lassoings has the unlinking number more than $r$ whereas $\mathrm{split}(L)=r$. 

\begin{figure}[!ht]
\begin{center}
\includegraphics[width=70mm]{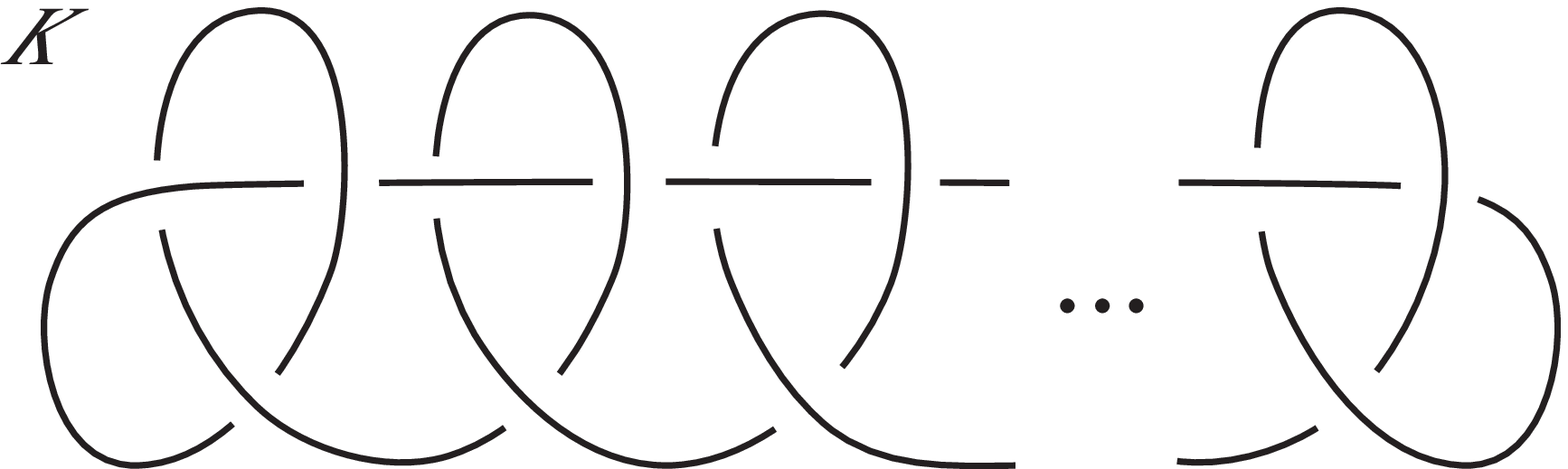}
\caption{}
\label{connected-tref}
\end{center}
\end{figure}
\end{eg}

\begin{ack}
The author is grateful to Reiko Shinjo for pointing out Adams' paper. 
She is also grateful to the members of Friday Seminar on Knot Theory in Osaka City University 
who assisted her in valuable comments and the tender encouragement. 
She especially thanks her advisor, Professor Akio Kawauchi for his help, support, and the valuable advices. 
She is partly supported by JSPS Research Fellowships for Young Scientists. 
\end{ack}

\maketitle

\end{document}